\documentclass{amsart}
\usepackage[utf8]{inputenc}
\usepackage[english]{babel}
\usepackage{verbatim}
\usepackage{graphicx}
\usepackage{amsmath,amsthm,amssymb}
\usepackage{afterpage} 
\usepackage{tcolorbox} 
\usepackage{float}
\usepackage{subfig}
\usepackage{enumerate}
\usepackage{todonotes}

\newtheorem{Theorem}{Theorem}
\newtheorem{Corollary}[Theorem]{Corollary}
\newtheorem{Proposition}[Theorem]{Proposition}
\newtheorem{Lemma}[Theorem]{Lemma}
\newtheorem{Definition}[Theorem]{Definition}
\newtheorem{Remark}[Theorem]{Remark}
\newtheorem{Property}[Theorem]{Property}

\newtheorem{Example}[Theorem]{Example}

\begin{document}

\title{Group invariant variational principles}

\author[J. Falc\'{o}]{Javier Falc\'{o}}
\address[Javier Falc\'{o}]{Departamento de An\'{a}lisis Matem\'{a}tico,
	Universidad de Valencia, Doctor Moliner 50, 46100 Burjasot (Valencia), Spain} \email{francisco.j.falco@uv.es}
\author[D. Isert]{Daniel Isert}
\address[Daniel Isert]{Departamento de An\'{a}lisis Matem\'{a}tico,
	Universidad de Valencia, Doctor Moliner 50, 46100 Burjasot (Valencia), Spain}
\email{isada@alumni.uv.es}

\thanks{The first author was supported by grant PID2021-122126NB-C33 funded by MCIN/AEI/10.13039/501100011033 and by “ERDF A way of making Europe”.}



\maketitle

\begin{abstract}
    In this paper we introduce a group invariant version of the well-known Ekeland variational principle. To achieve this, we define the concept of convexity with respect to a group and establish a version of the theorem within this framework. Additionally, we present several consequences of the group invariant Ekeland variational principle, including Palais-Smale minimizing sequences, the Brønsted-Rockafellar theorem, and a characterization of the linear and continuous group invariant functionals space. Moreover, we provide an alternative proof of the Bishop-Phelps theorem and proofs for the group-invariant Hahn-Banach separating theorems. Finally, we discuss some implications and applications of these results.

\end{abstract}
\tableofcontents

\section{Introduction}
Group invariant mappings have been studied through the years on the late century with, special emphasis in the symmetric multilinear mappings between Banach spaces. However, group invariant mappings in a general framework of infinite dimensional analysis, to our knowledge, classically have only been studied in the measure theory context, as we can see in \cite{RF}. Recently this general point of view has been thoroughly studied in some other fields of analysis, different from the measure theory, such as multivariable complex analysis or norm-attaining operators, for the latter ones see \cite{ArFaMa, ArFaGaMa, DFJ}.

The purpose of this paper is to give a group invariant version of the well-known Ekeland variational principle, which surprisingly, does not extend naturally to group invariant functionals like the Bishop-Phelps-Bollobás theorem presented in \cite{Falco}. We had to define the notion of convexity with respect to a group in order to prove this version of the theorem. Furthermore, we will present some consequences of the Ekeland variational principle such as the Palais-Smale minimizing sequences theorem, the Br\o nsted-Rockafellar theorem, a criteria to know when the set of group invariant points of a Banach space is again a Banach space, and a complete description of the linear and continuous group invariant functionals space. Finally we will give an alternative proof to the one given in \cite{Falco} of the Bishop-Phelps theorem, and we will also provide proofs of the group-invariant Hahn-Banach separating theorems.

We will begin this paper by presenting the notation we will use and setting the main definitions. Then we will move to study the Hahn-Banach separation theorems, extending the results considered in \cite[Section 3]{DFJ}. To continue we will give the main result of this paper, which is the generalization of Ekeland's variational principle for the group invariant functionals. To conclude, we will present some consequences of this result.

\section{Preliminaries}

Let us start by fixing some notation. We will denote a normed space by $X$, and the space of linear and continuous mappings from $X$ to $X$ by $\mathcal{L}(X)$. Hereinafter, our topological groups, $G$, are going to be subsets of this space, i.e, $G \subseteq \mathcal{L}(X)$. Hence, $G$ is reserved to denote a topological group of invertible bounded linear mappings with the topology endowed with the relative topology of $\mathcal{L}(X)$, also $g$ is reserved to denote the elements of the group $G$.

Let us begin by presenting the fundamental definitions of this paper.
\begin{Definition}
Let $X$, $Y$ be two normed spaces:
\begin{enumerate}
\item A point $x \in X$ is $G$-invariant, or invariant under the action of $G$, if $g(x) = x$ for all $g \in G$.

\item A set $K \subset X$ is $G$-invariant if for every $g \in G$, $g(K) = K$.

\item A mapping $f \colon X \to Y$ is $G$-invariant if for every $x \in X$ and every $g \in G$ we have that
\[
f(g(x)) = f(x).
\]
\end{enumerate}
\end{Definition}

We will denote by $X_{G}$, $B_{X}$, $S_{X}$ the group invariant points of $X$, the open unit ball, and the unit sphere of the space $X$ respectively. We will always assume that the norm of the space $X$, $\Vert \cdot \Vert$, is $G$-invariant.

We denote the set of continuous linear functionals that are $G$-invariant by
\[
X^{*}_{G} = \left\{f \in X^{*} ~ | ~ f \hbox{ is } G\hbox{-invariant}\right\}.
\]
Now we are going to define the symmetrization\footnote{The word symmetrization refers to the well-known process of applying Definition \ref{Simetritzacio de punt i funcional} to the group of permutations $G = \Sigma_{n}$ acting on $\mathbb{R}^{n}$.} operator by using the Bochner integral and the Haar measure. Roughly speaking this operator takes one point and by averaging its orbit under the action of $G$ produces a $G$-invariant point. This definition is a generalization of the Reynolds operator defined in \cite[Chapter 7]{CLO}.
Hereinafter, the measure used is always going to be the Haar measure associated to the compact group, and the integral, the Bochner integral. For the elemental notions on Haar measure we recommend \cite[Chapter 8]{Cohn}.
\begin{Definition}\label{Simetritzacio de punt i funcional}
Let $X$ be a Banach space and $G \subseteq \mathcal{L}(X)$ be a compact topological group acting on $X$. If $x \in X$ we define the symmetrization point of $x$ with respect to $G$, or the $G$-symmetrization point, as
\[
\overline{x} = \int_{G}g(x)d\mu(g).
\] 
Similarly if we have a functional $f \colon X \to \mathbb{R}$ we define the symmetric functional of $f$ with respect to $G$, or the $G$-symmetric functional, to be
\[
\overline{f}(x) = \int_{G}F(g(x))d\mu(g).
\]
\end{Definition}

We now define the notion of convexity and linearity with respect to the group that will be needed for the rest of the manuscript.
\begin{Definition}
Let $X$ be a Banach space and $G \subseteq \mathcal{L}(X)$ be a compact topological group acting on $X$. Let $\varphi \colon X \to \mathbb{R}$ be a function. We say that $\varphi$ is convex with respect to $G$ if 
\[
\varphi\left(\int_{G}g(x) d\mu(g)\right) \leq \int_{G}\varphi(g(x))d\mu(g) \quad \forall \, x \in X.
\]
And we say that $\varphi$ is linear with respect to $G$ if
\[
\varphi\left(\int_{G}g(x) d\mu(g)\right) = \int_{G}\varphi(g(x))d\mu(g) \quad \forall \, x \in X.
\]
\end{Definition}

It is easy to see that if $\varphi$ is $G$-invariant and convex, or linear, with respect to $G$ then
\[
\varphi(\overline{x}) \leq \varphi(x), \quad \hbox{ and, } \quad \varphi(\overline{x}) = \varphi(x),
\]
respectively. Observe also that, for a function, being convex or linear with respect to the group is a weaker condition than being convex or linear.

In section \ref{Secció de conseqüències} we will need to make use of the differentiability of a mapping in order to give a full description of the dual group invariant set of a Banach space.

We recall that if $X$ is a Banach space, $U \subseteq X$ is a nonempty and open set, and $F \colon U \to \mathbb{R}$ is a mapping, then the Gâteaux variation of $F$ at $u_{0} \in U$ in the direction of $h$ is the limit
\[
\lim_{t \to 0}\frac{F(u_{0}+th) - F(u_{0})}{t}
\]
if it exists and is finite. We denote this limit by $\delta F(u_{0})(h)$. If the Gâteaux variation of $F$ is linear and continuous, then $F$ is said to be Gâteaux differentiable. 

We recall that $F$ is Fréchet differentiable at $u_{0}$ if there exists a linear and continuous mapping $L \colon X \to \mathbb{R}$ such that
\[
\lim_{h \to 0}\frac{F(u_{0} + h) - F(u_{0}) - L(h)}{\Vert h \Vert}.
\]

It is well known that if $L$ satisfies the definition of Fréchet differentiable, then $L$ is unique, and we will denote it by $F'(u_{0})$. We note that the $G$-invariance of the function is preserved by the Gâteaux and Fréchet differentiability.

\section{Group invariant Hahn-Banach separating theorems}\label{Teoremes de separació}

Our goal in this section is to obtain some versions of the Hahn-Banach separation theorems for group invariant mappings.

Let's start recalling some definitions that will be used through this section.

If $X$ is a Banach space, an affine hyperplane is a subset $H$ of $X$ of the form
\[
H = \left\{x \in X ~ | ~ f(x) = \alpha\right\},
\]
where $f$ is a nonzero linear functional and $\alpha \in \mathbb{R}$ is a fixed constant.

If $X$ is a Banach space, and $A, B \subseteq X$, we say that the hyperplane $H = \left\{x \in X ~~ | ~~ f(x) = \alpha\right\}$ separates $A$ and $B$ if there exists some $\alpha \in \mathbb{R}$ so that
\[
f(x) \leq \alpha \quad \forall \, x \in A \hspace{0.2cm} \hbox{ and } \hspace{0.2cm} f(x) \geq \alpha \quad \forall \, x \in B.
\]
Also, we say that $H$ strictly separates $A$ and $B$ if there exists some $\epsilon > 0$ and some $\alpha \in \mathbb{R}$ such that
\[
f(x) \leq \alpha - \epsilon \quad \forall \, x \in A \hspace{0.2cm} \hbox{ and } \hspace{0.2cm} f(x) \geq \alpha + \epsilon \quad \forall \, x \in B.
\]

If $X$ is a Banach space and $C \subseteq X$ is an open convex set with $0 \in C$, for every $x \in X$ the functional
\[
p(x) = \inf\left\{\alpha > 0 ~ | ~ \alpha^{-1}x \in C\right\}
\]
is called the Minkowski functional of $C$.

We will make use of the following well known properties of the Minkowski functional.
\begin{Property}
The Minkowski functional of an open convex set $C$, $p$, satisfies the following properties
\begin{enumerate}[(i)]
    \item $p(\lambda x) = \lambda p(x)$ for every $x \in X$, $\lambda > 0$.
    
    \item $p(x + y) \leq p(x) + p(y)$ for every $x, y \in X$.
    
    \item There is a constant $M > 0$ such that $0 \leq p(x) \leq M\Vert x \Vert$ for every $x \in X$.
    
    \item $C = \left\{x \in X ~ | ~ p(x) < 1\right\}$.
\end{enumerate}
\end{Property}

The proof of the previous result can be found in \cite[Lemma 1.2]{Brezis}.

Moreover, in \cite[Theorem 3.1]{DFJ} it is proved that, if $C$ is a $G$-invariant subset, then the Minkowski functional is a $G$-invariant functional. Let us start with the following lemma.
\begin{Lemma}\label{Lema de la primera forma geometrica}
Let $X$ be a Banach space, $G \subseteq \mathcal{L}(X)$ be a compact topological group of isometries acting on $X$, and $C \subseteq X$ a $G$-invariant open convex set. If there exists a $G$-invariant point $x_{0} \in X\backslash C$, then we can find $f \in X^{*}_{G}$ such that $f(x) < f(x_{0})$ for all $x \in C$.
\end{Lemma}
\begin{proof}
Given $x \in C$ we know by \cite[Proposition 2.2]{DFJ} that $\overline{x} \in C$. Consider $C - \overline{x}$, which is a $G$-invariant set, and the point $x_{0} - \overline{x}$. We can suppose that $0 \in C$, otherwise we can do a translation. Note that, since $0 \in C$, we have that $x_{0} \neq 0$. Consider the subspace $H = \mathbb{R}x_{0}$ and the functional
\[
\begin{array}{cccl}
    h\colon & H & \to & \mathbb{R} \\
     & tx_{0} & \mapsto & t.
\end{array}
\]
Observe that $H$ is a $G$-invariant subspace, since $x_{0}$ is $G$-invariant, and so is the functional $h$ since
\[
h(g(tx_{0})) = h(tg(x_{0})) = h(tx_{0}) \quad \forall \, g \in G, \, \forall \, t \in \mathbb{R}.
\]
Moreover, by doing a study of cases when $t > 0$ and $t \leq 0$, is clear that $h(x) \leq p(x)$ for all $x \in H$, where $p$ is the Minkowski functional associated to $C$. Then, we can apply the $G$-invariant Hahn-Banach theorem, see \cite[Proposition 1]{Falco}, to obtain a functional $f \in X^{*}_{G}$ such that
\begin{enumerate}[(i)]
    \item $f_{|_{H}} = h$,
    
    \item $f(x) \leq p(x)$ for every $x \in X$.
\end{enumerate}
In particular $f(x_{0}) = 1$ and $f(x) \leq p(x) < 1$ for all $x \in C$.
\end{proof}

\begin{Theorem}\label{First Hahn-Banach geometric form}
Let $X$ be a Banach space and $G \subseteq \mathcal{L}(X)$ be a compact topological group os isometries acting on $X$. If $A$, $B$ are two convex and $G$-invariant sets, with $A$ open and $A \cap B = \emptyset$, then there exists a $G$-invariant hyperplane that separates $A$ and $B$.
\end{Theorem}
\begin{proof}
Define $C = A - B$ and notice that $C$ is convex and $G$-invariant, since $A$ and $B$ are $G$-invariant. Also, $C$ is open since we can write
\[
C = \bigcup_{y \in B} A - y,
\]
and $A$ is open. Moreover, $0$ is clearly $G$-invariant and does not belongs to $C$ since $A \cap B = \emptyset$. Therefore, we can apply Lemma \ref{Lema de la primera forma geometrica} to obtain a functional $f \in X^{*}_{G}$ with
\[
f(z) < 0 \quad \forall \, z \in C.
\]
Since $C = A - B$ and $f$ is linear, we deduce that
\[
f(x) < f(y) \quad \forall \, x \in A, \hspace{0.1cm} \forall \, y \in B.
\]
Then, there exists an $\alpha \in \mathbb{R}$ such that
\[
\sup_{x \in A} f(x) \leq \alpha \leq \inf_{y \in B}f(y).
\]
\end{proof}

The classical version of this theorem is usually called the first Hahn-Banach geometric form. Now we want to give a proof of the so-called second Hahn-Banach geometric form for $G$-invariant functionals. To do so we will require  to use the following lemma whose proof can be found in \cite[Theorem 3.1]{DFJ}.
\begin{Lemma}\label{Lema de la segona forma geometrica}
Let $X$ be a Banach space and $G \subseteq \mathcal{L}(X)$ be a compact topological group of isometries acting on $X$. If $C$ is a closed convex $G$-invariant subset of $X$ and $x_{0} \in X \backslash C$ is a $G$-invariant point, then there exists $f \in X^{*}_{G}$ such that
\[
f(x_{0}) > \sup_{x \in C}f(x).
\]
\end{Lemma}

The main difference with Lemma \ref{Lema de la primera forma geometrica} is the closeness of the set $C$, that makes the separation between the point and the set strict.
\begin{Theorem}
Let $X$ be a Banach space and $G \subseteq \mathcal{L}(X)$ be a compact topological group of isometries acting on $X$. If $A, B \subseteq X$ are nonempty, convex and $G$-invariant sets such that $A$ is closed, $B$ is compact and $A \cap B = \emptyset$, then there exists a $G$-invariant hyperplane that strictly separates $A$ and $B$.
\end{Theorem}
\begin{proof}
Define $C = A - B$ so that $C$ is convex, closed, $G$-invariant and $0 \notin C$. Hence, there exists some $r > 0$ sufficiently small such that $B(0, r) \cap C = \emptyset$.

Observe that by definition $B(0, r) = \left\{x \in X ~~ | ~~ \Vert x \Vert < r\right\}$, so since the norm is $G$-invariant, we can deduce the $G$-invariance of $B(0, r)$. Applying now Theorem \ref{First Hahn-Banach geometric form} we obtain that there exists $f \in X^{*}_{G}\backslash \left\{0\right\}$ such that
\[
f(x - y) \leq f(rz) \quad \forall \, x \in A, \, y \in B, \, \forall z \in B_{X}.
\]
By linearity of $f$, and symmetry of $B_{X}$, i.e, if $z \in B_{X}$ then $-z \in B_{X}$, we have that
\[
f(x - y) \leq -rf(z) \quad \forall \, x \in A, \, y \in B, \, \forall \, z \in B_{X},
\]
so
\[
f(x - y) \leq -r\Vert f \Vert \quad \forall \, x \in A, \, y \in B.
\]
Letting $\epsilon = \frac{1}{2}r\Vert f \Vert > 0$, we obtain
\[
f(x) + \epsilon < f(y) - \epsilon \quad \forall \, x \in A, \, y \in B.
\]
Now, we only have to choose $\alpha \in \mathbb{R}$ such that
\[
\sup_{x \in A}f(x) + \epsilon \leq \alpha \leq \inf_{y \in B}f(y) - \epsilon,
\]
and the results follows.
\end{proof}

To finish this section we will give a direct consequence of the previous result.
\begin{Corollary}
Let $X\ne\{0\}$ be a Banach space and $G \subseteq \mathcal{L}(X)$ be a compact topological group of isometries acting on $X$. Let $H \subseteq X$ be a $G$-invariant linear subspace such that $\overline{H} \neq X$. Then there exists a nonzero functional $f \in X^{*}_{G}$ such that
\[
\langle f, x \rangle = 0 \quad \forall \, x \in H.
\]
\end{Corollary}
\begin{proof}
Note that if $H=\{0\}$, since $X\neq\{0\}$, any non-zero functional of $X^{*}_{G}$ satisfies the conclusion.

Otherwise, observe that $\overline{H}$ is again $G$-invariant because so is $H$ and the points that are on the adherence are obtained as limit points of $G$-invariant sequences. Now let $x_{0} \in X \backslash \overline{H}$ be $G$-invariant. By the previous result we know that there exists $f \in X^{*}_{G}$, $f \not\equiv 0$, that strictly separates $\overline{H}$ and $\left\{x_{0}\right\}$, i.e.:
\[
\langle f, x \rangle < \alpha < \langle f, x_{0} \rangle \quad \forall \, x \in H,
\]
for some $\alpha \in \mathbb{R}$. 

Therefore, for any $x\in H, x\neq 0$ and using the last inequality we deduce that
\[
\langle f, x \rangle = 0 \quad \forall \, x \in H,
\]
since $\langle f, \lambda x \rangle = \lambda \langle f, x \rangle < \alpha$ for all $\lambda \in \mathbb{R}$.
\end{proof}

\section{Ekeland's variational principle for $G$-invariant functions}

In this section we are going to focus on the generalization of the Ekeland's variational principle to the context of $G$-invariant functionals. In the proof of this result it is required the following group invariant Cantor set theorem, whose proof can be obtained directly from the Cantor set theorem, and a small observation.

\begin{Proposition}
Let $(X, d)$ be a complete metric space, and let $G \subseteq \mathcal{L}(X)$ be a topological group acting on $X$. If $\left\{C_{n}\right\}_{n=1}^{+\infty}$ is a countable family of compact sets such that $C_{n} \neq \emptyset$, $C_{n+1} \subseteq C_{n}$ and $C_{n}$ is $G$-invariant for every $n \in \mathbb{N}$, then
\[
\bigcap_{n=1}^{+\infty}C_{n} \quad \hbox{is nonempty and } G \hbox{-invariant.}
\]
\end{Proposition}
\begin{proof}
By the Cantor set theorem we know that $\cap_{n=1}^{+\infty}C_{n} \neq \emptyset$. It only remains to show that $\cap_{n=1}^{+\infty}C_{n}$ is $G$-invariant, i.e.,
\[
g\left(\bigcap_{n=1}^{+\infty}C_{n}\right) = \bigcap_{n=1}^{+\infty}C_{n} \quad \forall \, g \in G.
\]
Let $x \in \cap_{n=1}^{+\infty}C_{n}$. Then, $x \in C_{n}$ for every $n \in \mathbb{N}$, but since every $C_{n}$ is $G$-invariant we get that
\[
g(x) \in C_{n} \quad \forall \, g \in G, \, \forall \, x \in C_{n}, \, \forall \, n \in \mathbb{N}.
\]
In particular, $g(C_{n}) \subseteq C_{n}$ for all $g \in G$, $n \in \mathbb{N}$. Therefore, $g\left( \bigcap_{n=1}^{+\infty}C_n\right)\subseteq C_n$ for all $g \in G$, $n \in \mathbb{N}$. Thus,
\[
g\left( \bigcap_{n=1}^{+\infty}C_n\right)\subseteq \bigcap_{n=1}^{+\infty}C_n \quad \forall g \in G.
\]

Using $g^{-1}$ we obtain the reverse inclusion, hence $\cap_{n=1}^{+\infty}C_{n}$ is $G$-invariant.
\end{proof}

As a direct consequence, we obtain that if the diameter of the sets converges to zero, then the intersection is a $G$-invariant point.

\begin{Corollary}\label{Intervals encaixats de Cantor}
Let $(X, d)$ be a complete metric space, and let $G \subseteq \mathcal{L}(X)$ be a topological group acting on $X$. If $\left\{C_{n}\right\}_{n=1}^{+\infty}$ is a countable family of compact sets such that $C_{n} \neq \emptyset$, $C_{n+1} \subseteq C_{n}$, $\delta(C_{n}) \to 0$ and $C_{n}$ is $G$-invariant for every $n \in \mathbb{N}$, then
\[
\bigcap_{n=1}^{+\infty}C_{n} = \left\{x\right\} \hbox{ is } G \hbox{-invariant,}
\]
where  $\delta(C_{n})$ denotes the diameter of the set $C_{n}$.
\end{Corollary}

Before we present the Ekeland variational principle, let us recall some required definitions.

If $X$ is a topological space, we say that $f \colon X \to \mathbb{R}$ is lower semicontinuous if for every $\lambda \in \mathbb{R}$ the set
\[
\left\{x \in X ~ | ~ f(x) \leq \lambda \right\}
\]
is closed.

We say that $f \colon X \to \mathbb{R}\cup\left\{-\infty, +\infty\right\}$ is proper if it has nonempty domain, never takes on the value $-\infty$, and it is not identically equal to $+\infty$. 

To continue we present the main result of this section, that is, the group invariant Ekeland's variational principle.

\begin{Theorem}\label{Principi variacional de Ekeland}
Let $X$ be a Banach space and $G \subseteq \mathcal{L}(X)$ be a compact topological group of isometries acting on $X$. Let $\varphi \colon X \to ]0, +\infty]$ be proper, lower semicontinuous, bounded below, $G$-invariant, and convex with respect to the group $G$. Then, given $\epsilon > 0$ and $\delta > 0$, there exists a $G$-invariant point $\Tilde{x} \in X$ such that
\[
\varphi(\Tilde{x}) < \varphi(x) + \epsilon\Vert x - \Tilde{x} \Vert \quad \forall \, x \in X, \hspace{0.1cm} x \neq \Tilde{x}.
\]
Moreover, if $x_{0} \in X$ satisfies that $\varphi(x_{0}) < \inf\left\{\varphi(x) ~ | ~ x \in X\right\} + \delta$, then we can choose $\Tilde{x}$ to be such that
\[
\Vert \overline{x_{0}} - \Tilde{x} \Vert < \frac{\delta}{\epsilon}.
\]
\end{Theorem}
\begin{proof}
Fix $\delta > 0$, we want to define a sequence $\left\{x_{n}\right\}_{n=1}^{+\infty} \subseteq X$ of $G$-invariant points that is convergent to our desired $\Tilde{x}$.

Define $b := \inf\left\{\varphi(x) ~ | ~ x \in X\right\}$, and let us verify that 
\[
b=\inf\left\{\varphi(x) ~ | ~ x \in X\right\} = \inf\left\{\varphi(x) ~ | ~ x \in X \hbox{ is } G\text{-invariant}\right\}.
\]
Since the set on the right hand side of the equality is smaller than the set on the left hand side, it is obvious that
\[
\inf\left\{\varphi(x) ~ | ~ x \in X\right\} \leq \inf\left\{\varphi(x) ~ | ~ x \in X \hbox{ is } G\text{-invariant}\right\}.
\]
For the other inequality, choose $x \in X$ and consider $\overline{x}$, the $G$-symmetric point of $x$.
Let's prove that $\varphi(\overline{x}) \leq \varphi(x)$. Using the definition of $\overline{x}$, the convexity of $\varphi$ with respect to the group and that $\varphi$ is $G$-invariant we have that
\[
\varphi(\overline{x}) = \varphi\left(\int_{G}g(x)d\mu(g)\right) \leq \int_{G}\varphi(g(x))d\mu(g) = \int_{G}\varphi(x)d\mu(g) = \varphi(x).
\]
Therefore, we obtain the desired equality, that is,
\[
\inf\left\{\varphi(x) ~ | ~ x \in X\right\} = \inf\left\{\varphi(x) ~ | ~ x \in X \hbox{ is } G\text{-invariant}\right\}.
\]
Since $b$ is an infimum, we can choose a $G$-invariant point $x_{0} \in X$ such that
\[
\varphi(x_{0}) < b + \frac{\delta}{2}.
\]
Inductively define, for $n = 0, 1, \ldots$,
\[
b_{n} = \inf\left\{t \in \mathbb{R} ~ | ~ \exists \, x \in X ~ G\text{-invariant} \colon (x, t) \in \hbox{epi }\varphi\cap\left((x_{n}, \varphi(x_{n})) + K_{\epsilon}\right)\right\},
\]
and choose $x_{n+1} \in X$, $G$-invariant, satisfying the following properties
\begin{equation}\label{Equació 3.1}
    (x_{n+1}, \varphi(x_{n+1})) \in (x_{n}, \varphi(x_{n})) + K_{\epsilon},
\end{equation}
\begin{equation}\label{Equació 3.2}
    \varphi(x_{n+1}) < b_{n} + \frac{\delta}{2^{n+2}},
\end{equation}
where we define
\[
K_{\epsilon} = \left\{(x, t) \in X \times \mathbb{R} ~ | ~ t \leq -\epsilon\Vert x \Vert\right\}.
\]
Observe that, since $\Vert \cdot \Vert$ is $G$-invariant, then so is $K_{\epsilon}$. Moreover by \eqref{Equació 3.1} we get: 
\[
(x_{n+1} - x_{n}, \varphi(x_{n+1}) - \varphi(x_{n})) \in K_{\epsilon}.
\]
By definition of $K_{\epsilon}$ we know that 
\[
\varphi(x_{n+1}) - \varphi(x_{n}) \leq -\epsilon \Vert x_{n+1} - x_{n} \Vert,
\]
so
\[
\varphi(x_{n}) - \varphi(x_{n + 1}) \geq \epsilon \Vert x_{n+1} - x_{n} \Vert \geq 0,
\]
and then
\[
\varphi(x_{n}) \geq \varphi(x_{n + 1}) \quad \forall n \in \mathbb{N}.
\]
Which means that $\left\{\varphi(x_{n})\right\}_{n=0}^{+\infty}$ is non increasing. 

Note that $\left\{(x_{n}, \varphi(x_{n})) + K_{\epsilon}\right\}_{n=0}^{+\infty}$ is a sequence of nested sets and $b_{n}$ is nondecreasing.

Now we will show that the sequence $\left\{x_{n}\right\}_{n=0}^{+\infty}$ is convergent to some point. By \eqref{Equació 3.1} and \eqref{Equació 3.2}, we have that
\[
b_{n-1} \leq b_{n} \leq \varphi(x_{n+1}) \leq \varphi(x_{n}) < b_{n-1} + \frac{\delta}{2^{n+1}} \quad \forall \, n \in \mathbb{N}.
\]
Therefore
\[
\varphi(x_{n}) - \varphi(x_{n-1}) < \frac{\delta}{2^{n+1}} \quad \forall \, n \in \mathbb{N}.
\]
And from this we can obtain that
\[
\epsilon\Vert x_{n} - x_{n+1} \Vert \leq \varphi(x_{n}) - \varphi(x_{n-1}) < \frac{\delta}{2^{n+1}} \quad \forall \, n \in \mathbb{N},
\]
hence
\[
\Vert x_{n} - x_{n+1} \Vert \leq \frac{\delta}{\epsilon 2^{n+1}} \quad \forall \, n \in \mathbb{N}.
\]
Thus, we obtain that $\left\{x_{n}\right\}_{n=0}^{+\infty}$ is a Cauchy sequence, but since $X$ is a Banach space, the sequence is convergent to some $\Tilde{x}\in X$.

Finally, we are going to check that we can apply Corollary \ref{Intervals encaixats de Cantor}. We already know that $(x_{n}, \varphi(x_{n})) + K_{\epsilon}$ are nested, and since $\varphi$ is lower semicontinuous, it is clear that $\hbox{epi }\varphi \cap ((x_{n}, \varphi(x_{n})) + K_{\epsilon})$ is closed for every $n\in \mathbb N$. Let us estimate their diameter in $(X \times \mathbb{R}, \rho)$ where $\rho(x, t) = \Vert x \Vert + |t|$ is the metric defined in the product space. Fix $n \in \mathbb{N} \backslash \left\{0\right\}$, given $(x, t) \in (x_{n}, \varphi(x_{n})) + K_{\epsilon}$, by the definition of $K_{\epsilon}$, we have that
\[
\epsilon\Vert x - x_{n} \Vert \leq \varphi(x_{n}) - t,
\]
and by the construction of the sequence
\[
\varphi(x_{n}) - t \leq \varphi(x_{n}) - b_{n} < \frac{\delta}{2^{n+1}}.
\]
Thus
\[
\Vert x - x_{n} \Vert \leq \frac{\delta}{\epsilon 2^{n+1}}.
\]
This proves that $\hbox{diam}\left(\hbox{epi }\varphi \cap ((x_{n}, \varphi(x_{n}))+ K_{\epsilon})\right) \to 0$. Applying then Corollary \ref{Intervals encaixats de Cantor} we get that
\[
\hbox{epi }\varphi \bigcap \left(\bigcap_{n=1}^{+\infty}(x_{n}, \varphi(x_{n})) + K_{\epsilon}\right) = \left\{(\Tilde{x}, \varphi(\Tilde{x}))\right\} \hbox{ and is } (G, Id)\hbox{-invariant}.
\]
In particular we have that $\Tilde{x}$ is $G$-invariant. 

By construction we have that $\hbox{epi }\varphi \cap ((\Tilde{x}, \varphi(\Tilde{x})) + K_{\epsilon}) = (\Tilde{x}, \varphi(\Tilde{x}))$, and in particular
\[
\varphi(\Tilde{x}) < \varphi(x) + \epsilon\Vert x - \Tilde{x} \Vert \quad \forall \, x \in X, \hspace{0.1cm} x \neq \Tilde{x}.
\]
Moreover, we have that 
\[
\Vert x_{0} - x_{n} \Vert \leq \Vert x_{0} - x_{1} \Vert + \dots + \Vert x_{n-1} - x_{n} \Vert \leq \frac{\delta}{\epsilon}\left(\frac{1}{2} + \dots + \frac{1}{2^{n}}\right),
\]
therefore $\Vert x_{0} - \Tilde{x} \Vert \leq \frac{\delta}{\epsilon}$.
\end{proof}

We are going now to give an example that shows that the convexity of the functional with respect to the group cannot be completely removed.
\begin{Example}
Consider the group $G = \left\{Id, -Id\right\}$ acting on $\mathbb{R}$. Then, the function 
\[
f(x) =
\begin{cases}
2x+1 & \hbox{if } -\frac{1}{2} \leq x \leq 0,\\
-2x+1 & \hbox{if } 0 \leq x \leq \frac{1}{2},\\
0 & \hbox{otherwise},
\end{cases}
\]
is not convex with respect to the group $G$. Observe that the only $G$-invariant point in $\mathbb{R}$ is $0$. Then, for every $\epsilon < 1$ and for $x = \frac{1}{2}$, we have that
\[
1 = f(0) \geq f\left(\frac{1}{2}\right) + \epsilon\frac{1}{2}.
\]
And thus, the conclusions of Theorem \ref{Principi variacional de Ekeland} does not hold.
\end{Example}

\section{Consequences of the Ekeland variational principle}\label{Secció de conseqüències}
\subsection{Palais-Smale minimizing sequences}
The first direct consequence of the Ekeland variational principle is the Palais-Smale minimizing sequences. Let us continue by generalizing this result.
\begin{Corollary}\label{Palais-Smale minimizing sequences}
Let $X$ be a Banach space and $G \subseteq \mathcal{L}(X)$ be a compact topological group of isometries. Let $\varphi \colon X \to \mathbb{R}$ be Gâteaux differentiable, bounded below, $G$-invariant and convex with respect to the group. Then, there exists a sequence $\left\{x_{n}\right\}_{n=1}^{+\infty} \subseteq X$ such that
\begin{enumerate}
    \item $x_{n}$ is $G$-invariant for all $n \in \mathbb{N}$,
    \item $\varphi(x_{n}) \to \inf\left\{\varphi(x) ~ | ~ x \in X\right\},$
    \item $\Vert \delta\varphi(x_{n})\Vert \to 0$.
\end{enumerate}
\end{Corollary}
\begin{proof}
Choose $\epsilon = \delta = \frac{1}{n}$. By Theorem \ref{Principi variacional de Ekeland} we can find a sequence $\left\{x_{n}\right\}_{n=1}^{+\infty} \subseteq X$ of $G$-invariant points such that
\[
\varphi(x_{n}) < \inf_{x \in X}\varphi(x) + \frac{1}{n},
\]
and such that
\[
\varphi(x_{n}) \leq \varphi(x) + \frac{1}{n}\Vert x - x_{n} \Vert \quad \forall \, x \in X.
\]
Therefore
\[
\varphi(x_{n}) - \varphi(x) \leq \frac{1}{n}\Vert x - x_{n} \Vert \quad \forall \, x \in X.
\]
Now let's check that $\Vert \delta\varphi(x_{n}) \Vert \to 0$. We know that
\[
\Vert \delta\varphi(x_{n}) \Vert = \sup_{y \in X, \Vert y \Vert \leq 1}|\delta\varphi(x_{n})(y)|.
\]
Fix $y \in B_{X}$. Defining $x = x_{n} + ty$, for $t > 0$, we have that
\[
\varphi(x_{n}) - \varphi(x_{n} + ty) < \frac{1}{n}\Vert ty \Vert,
\]
therefore
\[
\frac{\varphi(x_{n}) - \varphi(x_{n} + ty)}{t} \leq \frac{1}{n}\Vert y \Vert \quad \forall \, y \in B_{X}.
\]
Taking now limits when $t \to 0$ we get that
\[
-\delta\varphi(x_{n})(y) = \lim_{t \to 0}\frac{\varphi(x_{n}) - \varphi(x_{n} + ty)}{t} \leq \frac{\Vert y \Vert}{n} \quad \forall \, y \in B_{X}.
\]
Since we have the inequality for every $y$, considering $-y$ and using the linearity of the Gâteaux differential, we have obtained that
\[
|\delta\varphi(x_{n})(y)| \leq \frac{\Vert y \Vert}{n} \leq \frac{1}{n} \quad \forall \, y \in B_{X}.
\]
Taking now limits when $n \to +\infty$ we obtain the desired result.
\end{proof}

\begin{Corollary}
Let $X$ be a Banach space and $G \subseteq \mathcal{L}(X)$ be a compact topological group of isometries. Let $f \colon X \to \mathbb{R}\cup\left\{+\infty\right\}$ be a lower semicontinuous function, Gâteaux differentiable, bounded below, $G$-invariant, convex with respect to $G$, and so that there exists constants $k, c > 0$ with
\[
f(x) \geq k\Vert x \Vert + c \quad \forall \, x \in X.
\]
Then, the range of $\delta f(x)$ is dense in $kB_{G}$, where $B_{G}$ is the closed unit ball in $X_{G}^{*}$.
\end{Corollary}
\begin{proof}
Choose $u \in kB_{G}$. We want to see if there exists a sequence of $G$-invariant points $\left\{u_{n}\right\}_{n=1}^{+\infty} \subseteq X$ such that
\[
\lim_{n \to +\infty}\Vert \delta f(u_{n}) - u \Vert = 0.
\]
Define
\[
\begin{array}{cccl}
    h \colon & X & \to & \mathbb{R} \\
     & x & \mapsto & f(x) - \langle u, x \rangle.
\end{array}
\]
It is easy to see that $h$ is Gâteaux differentiable, $G$-invariant and convex with respect to $G$. Note that
\begin{align*}
\inf_{x \in X} h(x) & = \inf_{x \in X} f(x) - \langle u, x \rangle \\
& \geq \inf_{x \in X} k\Vert x \Vert + c - \langle u, x \rangle \\
& \geq \inf_{x \in X} k\Vert x \Vert + c - \Vert u \Vert\hspace{0.1cm}\Vert x \Vert \\
& = \inf_{x \in X} \left( k - \Vert u \Vert\right)\Vert x \Vert + c \\
& \geq c > -\infty,
\end{align*}
where we have used the hypothesis on the first inequality, and that $u \in kB_{G}$ on the last line. Then, $h$ is bounded below, and it satisfies the hypothesis of Corollary \ref{Palais-Smale minimizing sequences}, therefore there exists a sequence of $G$-invariant points $\left\{u_{n}\right\}_{n=1}^{+\infty}$ such that $\lim_{n \to +\infty}\Vert \delta h(u_{n}) \Vert = 0$, in particular
\[
\lim_{n \to +\infty}\Vert \delta f(u_{n}) - u \Vert = 0.
\]
\end{proof}

\begin{Corollary}
Let $X$ be a Banach space and $G \subseteq \mathcal{L}(X)$ be a compact topological group of isometries. Let $f \colon X \to \mathbb{R}\cup\left\{+\infty\right\}$ be a lower semicontinuous function, Gâteaux differentiable, bounded below, $G$-invariant, convex with respect to $G$, and suppose that there exists $\phi \colon \mathbb{R}_{+} \to \mathbb{R}\cup\left\{+\infty\right\}$ such that $\frac{\phi(t)}{t} \to +\infty$ when $t \to +\infty$ and
\[
f(x) \geq \phi\left(\Vert x \Vert\right) \quad \forall x \in X.
\]
Then, the range of $\delta f(x)$ is dense in $X^{*}_{G}$.
\end{Corollary}
\begin{proof}
It is clear that for all $k > 0$, there exists $c \in \mathbb{R}$ such that $f$ satisfies that
\[
f(x) \geq k\Vert x \Vert + c.
\]
Then, $\delta f(x)$ is dense on every closed ball of $X^{*}_{G}$
\end{proof}

\subsection{Description of $X^{*}_{G}$}

Let us start giving a lemma that will be useful on the proof of Theorem \ref{Dual grup invariant}.
\begin{Lemma}\label{Lema previ al corolari del dual dels G-invariants}
Let $X$ be a Banach space, $G \subseteq \mathcal{L}(X)$ be a compact topological group acting on $X$, and $\varphi \colon X \to \mathbb{R}$ be $G$-invariant. If $\varphi$ is linear with respect to $G$, then $\frac{1}{\varphi^{2}}$ is convex with respect to $G$.
\end{Lemma}
\begin{proof}
Define $\psi \colon X \to \mathbb{R}$ as follows
\[
\psi(x) = 
\begin{cases}
\frac{1}{\varphi^{2}(x)} & \hbox{if } \varphi(x) \neq 0,\\
+\infty & \hbox{if } \varphi(x) = 0.
\end{cases}
\]
Observe that
\[
\psi\left(\int_{G}g(x)d\mu(g)\right) \leq \int_{G}\psi(g(x))d\mu(g).
\]
Which is equivalent, by definition, to
\[
\frac{1}{\varphi^{2}\left(\int_{G}g(x)d\mu(g)\right)} \leq \displaystyle{\int_{G}\frac{1}{\varphi^{2}(g(x))}d\mu(g)}.
\]
By $G$-invariance of $\varphi$ and definition of the symmetrized point, this is equivalent to
\[
\frac{1}{\varphi^{2}\left(\overline{x}\right)} \leq \displaystyle{\int_{G}\frac{1}{\varphi^{2}(x)}d\mu(g) = \frac{1}{\varphi^{2}(x)}}.
\]
Hence,
\[
\varphi^{2}(x) \leq \varphi^{2}(\overline{x}).
\]
By doing here a distinction of cases when $\varphi \geq 0$ and when $\varphi \leq 0$, we obtain the desired result.
\end{proof}

\begin{Theorem}\label{Dual grup invariant}
Let $X$ be a Banach space and $G \subseteq \mathcal{L}(X)$ be a compact topological group of isometries acting on $X$. Let $\varphi$ be a continuous function, Gâteaux differentiable, bump, $G$-invariant and linear with respect to $G$. Then
\[
X^{*}_{G} = \overline{\text{Span}}\left\{\partial\varphi(x) ~ | ~ x \in X_{G}\right\}.
\]
\end{Theorem}
\begin{proof}
Define $\psi \colon X \to \mathbb{R}\cup\left\{+\infty\right\}$ as follows
\[
\psi(x) = 
\begin{cases}
\frac{1}{\varphi(x)^2} & \hbox{if } \varphi(x) \neq 0,\\
+\infty & \hbox{if } \varphi(x) = 0.
\end{cases}
\]
Let $f \in X^{*}_{G}$. Then it is clear that $\psi(x) - f(x)$ is lower semicontinuous and bounded below function. Notice also that $\psi$ is Gâteaux differentiable and by Lemma \ref{Lema previ al corolari del dual dels G-invariants} is convex with respect to the group. Therefore, by Theorem \ref{Principi variacional de Ekeland}, we have that there exists a $G$-invariant point $x_{0} \in X$ such that for every $h \in X$ and $t > 0$ we have that
\[
\psi(x_{0} + th) - f(x_{0} + th) \geq \psi(x_{0}) - f(x_{0}) - \epsilon t \Vert h \Vert.
\]
Hence
\[
\frac{\psi(x_{0} + th) - \psi(x_{0})}{t} \geq \frac{f(x_{0} + th) - f(x_{0})}{t} - \epsilon\Vert h \Vert,
\]
and by linearity of $f$ we obtain
\[
\frac{\psi(x_{0} + th) - \psi(x_{0})}{t} \geq f(h) - \epsilon\Vert h \Vert.
\]
Taking now limits when $t \to 0^{+}$ we finally have that
\[
\partial\psi(x_{0})(h) = \lim_{t \to 0^{+}}\frac{\psi(x_{0} + th) - \psi(x_{0})}{t} \geq f(h) - \epsilon\Vert h \Vert \quad \forall \, h \in X.
\]
Considering $-h$ and using the linearity of $f$ and $\partial \psi$, it follows that
\[
|\partial\psi(x_{0})(h) - f(h)| \leq \epsilon\Vert h \Vert \quad \forall \, h \in X.
\]
Finally, we conclude that
\[
\left\Vert -2\frac{\partial\varphi(x_{0})}{\varphi(x_{0})^{3}} - f \right\Vert = \Vert \partial\psi(x_{0}) - f \Vert \leq \epsilon \quad \forall \, \epsilon > 0.
\]
\end{proof}

\subsection{A characterization of the completeness of $X_{G}$.}
Our next result is based on the characterization given by Sullivan on 1981 about the completeness of a metric space in \cite{Sullivan}, and tells us exactly when the $G$-invariant normed spaces are indeed Banach spaces.
\begin{Theorem}
Let $\left(X, \Vert \cdot \Vert\right)$ be a normed space, and $G \subseteq \mathcal{L}(X)$ be a compact topological group of isometries acting on $X$. $X_{G}$ is a Banach space if, and only if, for all $f \colon X \to \mathbb{R}$ that is bounded below, Lipschitz continuous, $G$-invariant, convex with respect to $G$ ,and for all $\epsilon > 0$, there exists $x_{0} \in X_{G}$ such that
\[
f(x_{0}) \leq \inf_{X}f + \epsilon,
\]
and
\[
f(x_{0}) \leq f(x) + \epsilon\Vert x - x_{0} \Vert \quad \forall \, x_{0} \neq x \in X.
\]
\end{Theorem}
\begin{proof}
The direct implication is clear by applying Theorem \ref{Principi variacional de Ekeland} with $\delta = \epsilon^{2}$.

For the sufficient condition, fix a Cauchy sequence $\left\{x_{n}\right\}_{n=1}^{+\infty} \subseteq X_{G}$. Then for all $\epsilon < 1$ there exists $N_{0} \in \mathbb{N}$ such that
\[
\Vert x_{n} - x_{m} \Vert < \epsilon \quad \forall n > m \geq N_{0}.
\]
Define now the functional
\[
\begin{array}{cccl}
    f \colon & X & \to & \mathbb{R} \\
     & x & \mapsto & \displaystyle{\hspace{-0.3cm}\lim_{n \to +\infty}\Vert x_{n} - \overline{x} \Vert}.
\end{array}
\]
Observe that it is non negative and satisfies that $\inf\left\{f(x) ~ | ~ x \in X\right\} = 0$. Let's see that it is well defined. We know that
\[
\left| \, \Vert x_{n} - \overline{x} \Vert - \Vert x_{m} - \overline{x} \Vert \, \right| \leq \Vert x_{n} - x_{m} \Vert < \epsilon \quad \forall \, n > m > N_{0},
\]
so $\left\{\Vert x_{n} - \overline{x} \Vert\right\}_{n=1}^{+\infty} \subseteq \mathbb{R}$ is a Cauchy sequence, hence it is convergent and then $f$ is well defined. Now we are going to see that $f$ is Lipschitz-continuous. Indeed, given $x, y \in X$ and $n \in \mathbb{N}$ we know that
\[
\left| \, \Vert x_{n} - \overline{x} \Vert - \Vert x_{n} - \overline{y} \Vert \, \right| \leq \Vert \overline{x} - \overline{y} \Vert = \left\Vert \int_{G}(g(x) - g(y))d\mu(g) \right\Vert \leq
\]\[
\leq \int_{G}\Vert g(x) - g(y)\Vert d\mu(g) = \int_{G}\Vert x - y\Vert d\mu(g) = \Vert x - y \Vert,
\]
where we have used here the definition of the $G$-symmetrization and the fact that the norm is $G$-invariant. Then taking limits when $n \to +\infty$, we deduce that
\[
\left|f(x) - f(y)\right| \leq \Vert x - y \Vert \quad \forall \, x, y \in X,
\]
i.e., $f$ is $1$-Lipschitz-continuous. Finally, let's see that $f$ is $G$-invariant and linear with respect to $G$. Pick $g \in G$, then
\[
f(g(x)) = \lim_{n \to +\infty}\Vert x_{n} - g(\overline{x}) \Vert = \lim_{n \to + \infty}\Vert x_{n} - \overline{x} \Vert = f(x),
\]
so $f$ is $G$-invariant. 

Notice that
\[
f(\overline{x}) = \lim_{n \to +\infty}\Vert x_{n} - \overline{\overline{x}} \Vert = \lim_{n \to +\infty}\Vert x_{n} - \overline{x} \Vert = f(x),
\]
where we have used that the sequence is made of $G$-invariant points, and the $G$-invariance of the norm. Moreover, since $\left\{x_{n}\right\}_{n=1}^{+\infty}$ is a Cauchy sequence, we have that
\[
f(x_{m}) = \lim_{n \to +\infty}\Vert x_{n} - x_{m} \Vert < \epsilon \quad \forall \, m \geq N_{0}.
\]
Hence,
\[
f\left(\int_{G}g(x)d\mu(g)\right) = f(\overline{x}) = f(x) = \int_{G}f(g(x))d\mu(g)
\]
Now, we can apply our hypothesis to obtain a $G$-invariant point, say $x_{0} \in X$, such that
\begin{enumerate}
    \item $f(x_{0}) \leq \epsilon$,
    
    \item $f(x_{0}) \leq f(x) + \epsilon\Vert x - x_{0} \Vert$ for all $x_{0} \neq x \in X$.
\end{enumerate}
By using $(2)$ with $x = x_{m}$, we obtain that
\begin{equation}\label{Completessa de l'espai}
f(x_{0}) \leq f(x_{m}) + \epsilon\Vert x_{m} - x_{0} \Vert = f(x_{m}) + \epsilon\Vert x_{m} - \overline{x_{0}} \Vert \quad \hbox{for each } m \in \mathbb{N}
\end{equation}
Notice now that
\[
f(x_{m}) \to 0 \hbox{ when } m \to +\infty,
\]
and
\[
\Vert x_{m} - \overline{x_{0}} \Vert \to \lim_{m \to +\infty}\Vert x_{m} - \overline{x_{0}} \Vert = f(x_{0}).
\]
So, taking limits over $m$ in \eqref{Completessa de l'espai} we have that
\[
f(x_{0}) \leq \epsilon f(x_{0}),
\]
which means that $(1 - \epsilon)f(x_{0}) \leq 0$. But from the start we know that $\epsilon < 1$ and that $\inf_{x \in X}f(x) = 0$, hence, $f(x_{0}) = 0$, which means that $\Vert x_{n} - \overline{x_{0}} \Vert \to 0$ when $n \to +\infty$. Finally, it is clear now that $x_{0}$ is $G$-invariant, and $x_{n} \to x_{0}$ when $n \to +\infty$. Thus, $X_{G}$ is a Bnach space.
\end{proof}

\begin{Remark}
    Note that the fact that $X_G$ is a Banach space does not guarantee that $X$ is a Banach space. Indeed, if we consider for instance the space $c_{00}$ consisting of all real sequences of finite support and the group of permutations of the first $n$ coordinates $G=\Sigma_n$ with $n$ an arbitrary but fixed natural number, it is clear that $G$ has a natural action on $X$ by 
\[\varphi(x_1,x_2,\ldots,x_n,x_{n+1},\ldots)=(x_{\varphi^{-1}(1)},x_{\varphi^{-1}(2)},\ldots,x_{\varphi^{-1}(n)},x_{n+1},\ldots),
\]
for every $\varphi \in G$ and every $(x_1,x_2,\ldots)\in c_{00}$. Then $c_{00}$ is not a Banach space, but $(c_{00})_G$ is clearly a Banach space.
\end{Remark}

\subsection{A proof of the Bishop-Phelps theorem by using Ekeland's variational principle}

Here we are going to present a proof of the group invariant Bishop-Phelps theorem, by using the group invariant Ekeland's variational principle. This alternative proof is based on the proof of the Bishop-Phelps theorem presented in \cite[Theorem 3.18]{Phelps}.
\begin{Theorem}
Let $X$ be a real Banach space and $G \subseteq \mathcal{L}(X)$ be a compact topological group of isometries acting on $X$. If $C \subseteq X$ is a convex, closed, bounded and $G$-invariant subset, then the $G$-invariant norm-attaining functionals in $C$ are dense in $X^{*}_{G}$.
\end{Theorem}
\begin{proof}
Choose $f \in X^{*}_{G}$, $\epsilon > 0$, and define $\Tilde{f} \colon X \to \mathbb{R}\cup\left\{+\infty\right\}$
\[
\Tilde{f}(x) = 
\begin{cases}
-f(x) & \hbox{if } x \in C,\\
+\infty & \hbox{if } x \in X \backslash C.
\end{cases}
\]
Notice that $\Tilde{f}$ is proper, lower semicontinuous (because $f$ is continuous), and bounded below. We have to show that $\Tilde{f}$ is $G$-invariant. Indeed, when $x \in C$, we recall that by the convexity of $C$, $\overline{x} \in C$, so we can apply the $G$-invariance of $f$ to obtain that
\[
\Tilde{f}(g(x)) = -f(g(x)) = -f(x) = \Tilde{f}(x) \quad \forall \, g \in G,
\]
and when $x \in X \backslash C$, in which case, by the $G$-invariance of $C$ we have that
\[
\Tilde{f}(x) = +\infty = \Tilde{f}(g(x)),
\]
since $C$ is convex. Let us continue by showing that $\Tilde{f}$ is convex with respect to $G$. Indeed, for $x \in C$, by linearity and $G$-invariance of $f$ we obtain that
\[
\Tilde{f}(\overline{x}) = -f(\overline{x}) = -f\left(\int_{G}g(x)d\mu(g)\right) = -f(x) = \Tilde{f}(x).
\]
Now, fix $x \in X \backslash C$. Then, $\overline{x} \in X \backslash C$, so
\[
\Tilde{f}(\overline x) = +\infty = \Tilde{f}(x).
\]
Hence, $\Tilde{f}$ is linear with respect to $G$, and in particular it is convex with respect to $G$. Applying now Theorem \ref{Principi variacional de Ekeland} we obtain a $G$-invariant point $x_{0} \in C$ such that
\[
\Tilde{f}(x_{0}) \leq \Tilde{f}(x) + \epsilon \Vert x - x_{0} \Vert \quad \forall \, x \in X.
\]
By definition of $\Tilde{f}$ we have that
\[
f(x_{0}) + \epsilon\Vert x - x_{0} \Vert \geq f(x) \quad \forall \, x \in C.
\]
Define the sets
\[
K_{1} = \left\{(x, t) ~~ | ~~ x \in C, \, t \leq f(x)\right\} \subseteq X \times \mathbb{R},
\]\[
K_{2} = \left\{(x, t) ~~ | ~~ t \geq f(x_{0}) + \epsilon \Vert x - x_{0} \Vert\right\},
\]
and observe that both of them are $(G, Id)$-invariant since the first one is the epigraph of the function $f$, and the second one is a $G$-invariant set by the $G$-invariance of the norm and $x_{0}$. Moreover, the interior of $K_{2}$ is nonempty and disjoint from $K_{1}$. Therefore, by Theorem \ref{First Hahn-Banach geometric form}, there exists a $(G, Id)$-invariant functional in $X^{*} \times \mathbb{R}$ that separates $K_{1}$ and $K_{2}$, i.e., there exists $h \in X^{*}\backslash\left\{0\right\}$ and $\beta \in \mathbb{R}$ such that
\[
\langle h, x \rangle + \alpha t \leq \beta \quad \hbox{if } (x, t) \in K_{1},
\]
and
\[
\langle h, x \rangle + \alpha t \geq \beta \quad \hbox{if } (x, t) \in K_{2},
\]
for some $\alpha \in \mathbb{R}$. Observe that $\alpha$ cannot be negative, otherwise $\langle h, x \rangle + \alpha t < \beta$ if $(x, t) \in K_{2}$ and $t$ is large enough. On the same way, $\alpha$ cannot be equal to $0$, otherwise $\langle h, x \rangle \geq \beta$ for all $x \in X$, meaning that $h = 0$. Then, $\alpha > 0$, and we can normalize it to be $\alpha = 1$.

Notice now that, since $(x_{0}, f(x_{0})) \in K_{1} \cap K_{2}$ we have that $\langle h, x_{0} \rangle + f(x_{0}) = \beta$. Therefore
\[
\langle h, x \rangle + f(x) \leq \langle h, x_{0} \rangle + f(x_{0}) \quad \forall \, x \in C.
\]
From where we deduce that $h + f$ attains its maximum on $C$ at $x_{0}$. On the other hand, if $x \in X$ and $t = f(x_{0}) + \epsilon \Vert x - x_{0} \Vert$, then
\[
\langle h, x_{0} \rangle + f(x_{0}) = \beta \leq \langle h, x \rangle + f(x_{0}) + \epsilon \Vert x - x_{0} \Vert \quad \forall \, x \in X,
\]
since $(x, t) \in K_{2}$.
Hence
\[
\langle h, x_{0} - x \rangle \leq \epsilon\Vert x - x_{0} \Vert,
\]
so
\[
|\langle h, z \rangle| \leq \epsilon\Vert z \Vert \quad \forall \, z \in X,
\]
and then $\Vert h \Vert \leq \epsilon$.

Defining now $\Tilde{h} = f + h$, it is clear that $\Tilde{h}$ attains its maximum on $C$, is $G$-invariant and
\[
\Vert \Tilde{h} - f \Vert \leq \epsilon.
\]
\end{proof}

\begin{Remark}
The previous result in the complex case holds if $C$ is the unit ball.
\end{Remark}

A closer look at the previous proof shows that, indeed, the following equivalence of the Bishop-Phelps theorem holds.
\begin{Corollary}
Let $X$ be a real Banach space, $G \subseteq \mathcal{L}(X)$ be a compact topological group of isometries acting on $X$, and $C$ a nonempty, closed, convex and $G$-invariant subset of $X$. For $\epsilon > 0$, suppose that $f \colon C \to \mathbb{R}$ is lower semicontinuous, convex, bounded below and $G$-invariant. Then, there exists $h \in X^{*}_{G}$ such that $\Vert h \Vert \leq \epsilon$, and $f + h$ attains its maximum at some $G$-invariant point $x_{0} \in C$.
\end{Corollary}

\subsection{Br\o nsted-Rockafellar theorem}

In this section we want to give a group invariant version of the Br\o nsted-Rockafellar theorem, see \cite[Theorem 3.17]{Phelps}. Before stating the result, let us recall the notions of $\epsilon$-subdifferential, subdifferential and subgradient.

Let $X$ be a Banach space and $f \colon X \to \mathbb{R}\cup\left\{+\infty\right\}$ a proper, convex and lower semicontinuous function. Given $x_{0} \in \text{Dom}(f)$ and $\epsilon > 0$, the $\epsilon$-subdifferential of $f$ at $x_{0}$ is defined as
\[
\partial_{\epsilon}f(x_{0}) = \left\{h \in X^{*} ~~ | ~~ \langle h, x - x_{0} \rangle \leq f(x) - f(x_{0}) + \epsilon \hspace{0.2cm} \forall \, x \in X \right\}.
\]
If $x_{0} \notin \text{Dom}(f)$ we say that $\partial f(x_{0}) = \emptyset$.

The subdifferential of $f$ at $x_{0}$ is defined as
\[
\partial f(x_{0}) = \left\{h \in X^{*} ~~ | ~~ \langle h, x - x_{0} \rangle \leq f(x) - f(x_{0}) \hspace{0.2cm} \forall \, x \in X \right\}.
\]
If $x_{0} \notin \text{Dom}(f)$ we say that $\partial f(x_{0}) = \emptyset$. Also we say that $h \in X^{*}$ is a subgradient of $f$ at $x_{0} \in \text{Dom}(f)$ if
\[
\langle h, x - x_{0} \rangle \leq f(x) - f(x_{0}) \quad \forall \, x \in X.
\]

We will also need to recall the definition of adjoint operator.

If $X, Y$ are two Banach spaces, and $T \colon X \to Y$ is a bounded operator, the adjoint operator of $T$, denoted by $T^{*} \colon Y^{*} \to X^{*}$, is the operator defined by
\[
\langle x, T^{*}y^{*} \rangle = \langle Tx, y^{*} \rangle \quad \forall \, x \in X\hbox{, } y^{*} \in Y^{*}.
\]

Let us start with some elemental properties of the adjoint operator.
\begin{Proposition}
If $G \subseteq L(X)$ is a compact topological group, then so is $G^{*} \subseteq L(X^{*})$.
\end{Proposition}

The proof of this result follows from the norm-to-norm continuity of the adjoint operator and the compacity of $G$.

\begin{Property}
Let $X$ be a Banach space, $G \subseteq \mathcal{L}(X)$ be a compact topological group acting on $X$, $f \colon X \to \mathbb{R}\cup\left\{+\infty\right\}$ be a proper function, and $x_{0} \in \text{Dom}(f)$ be a $G$-invariant point. Then, the $\epsilon$-subdifferential and the subdifferential of $f$ at $x_{0}$ are $G^{*}$-invariant.
\end{Property}
\begin{proof}
We are going to show the $G^{*}$-invariance of the $\partial f(x_{0})$. Let $g^{*} \in G^{*}$, we want to see that if we choose $g^{*} \in G^{*}$, then $g^{*}(\partial f(x_{0})) = \partial f(x_{0})$. Observe that, by the definition of the adjoint of an operator, the linearity of $g$ and the the $G$-invariance of $x_{0}$ we have that
\[
\langle g^{*}u, x - x_{0} \rangle = \langle u, g(x - x_{0}) \rangle = \langle u, g(x) - g(x_{0}) \rangle = \langle u, g(x) - x_{0} \rangle.
\]
And now, since $u \in \partial f(x_{0})$, and $f$ is $G$-invariant it is clear that
\[
\langle g^{*}u, x - x_{0} \rangle = \langle u, g(x) - x_{0} \rangle \leq f(g(x)) - f(x_{0}) = f(x) - f(x_{0}).
\]
This means that $g^{*}u \in \partial f(x_{0})$, so we deduce that $g^{*}(\partial f(x_{0})) \subseteq \partial f(x_{0})$. Using the inverse mapping of $g^{*}$, we deduce the other inclusion. Then, we have that $g^{*}(\partial f(x_{0})) = \partial f(x_{0})$ for every $g^{*} \in G^{*}$.
\end{proof}

\begin{Remark}
An operator $T \in L(X, Y)$ is $G$-invariant if, and only if, $T^{**}$ is $G$-invariant, see \cite[Proposition 2.4]{DFJ} for the details.
\end{Remark}

Now, we are going to give, without proof, two results that will prove crucial for the proof of the Br\o nsted-Rockafellar theorem. The first one is the so called Moreau-Rockafellar theorem, which says, roughly-speaking, that the sum of the subdifferentials is the subdifferential of the sum.  For the proof of this results one can see \cite[Theorem 3.16]{Phelps}.
\begin{Theorem}\label{Moreau-Rockafellar}
Let $X$ be a Banach space and let $f, g \colon X \to \mathbb{R}\cup\left\{+\infty\right\}$ be two proper, convex and lower semiconitnuous functions. Suppose that $\text{Dom}(f) \cap \text{Dom}(g) \neq \emptyset$, then
\[
\partial f(x) + \partial g(x) \subseteq \partial(f + g)(x) \quad \forall \, x \in \text{Dom}(f + g).
\]
Moreover, if there exists one point in $\text{Dom}(f) \cap \text{Dom}(g)$ for which one of the two functions is continuous, then we have that
\[
\partial f(x) + \partial g(x) = \partial(f + g)(x) \quad \forall \, x \in \text{Dom}(f + g).
\]
\end{Theorem}

And as an easy Corollary of this result, we have the following.
\begin{Corollary}\label{Corollary Moreau-Rockafellar}
Let $X$ be a Banach space, $f \colon X \to \mathbb{R}\cup\left\{+\infty\right\}$ be proper, convex and lower semiconitnuous, and $h \in X^{*}$. For $x_{0} \in \text{Dom}(f)$:
\[
\partial f(x_{0}) + h = \partial(f + h)(x_{0}).
\]
\end{Corollary}

With this brief introduction we can move now to see the $G$-invariant version of the Br\o nsted-Rockafellar theorem.
\begin{Theorem}
Let $X$ be a Banach space and $G \subseteq \mathcal{L}(X)$ be a compact topological group of isometries acting on $X$. Suppose that $f \colon X \to \mathbb{R}\cup\left\{+\infty\right\}$ is a proper, convex, lower semicontinuous and $G$-invariant function. Then for any $G$-invariant point $x_0 \in Dom(f)$ and any $G$-invariant functional $x_{0}^{*} \in \partial_{\epsilon}f(x_{0})$, and for all  $\epsilon, \lambda > 0$, there exists a $G$-invariant point $z \in \text{Dom}(f)$ and a functional $x^{*} \in X^{*}_{G}$ such that
\[
x^{*} \in \partial f(x), \hspace{0.3cm} \Vert z - x_{0} \Vert \leq \frac{\epsilon}{\lambda}, \hspace{0.3cm} \Vert x^{*} - x_{0}^{*} \Vert \leq \lambda.
\]
\end{Theorem}
\begin{proof}
Fix a $G$-invariant point $x_{0} \in \text{Dom}(f)$, $\epsilon > 0$ and $\lambda > 0$. Then, we know that $\partial_{\epsilon}f(x_{0})$ is nonempty and $G^{*}$-invariant. Consider a $G$-invariant functional $x_{0}^{*} \in \partial_{\epsilon}f(x_{0})$ and define
\[
\begin{array}{cccc}
    \phi \colon & X & \to & \mathbb{R}\cup\left\{+\infty\right\} \\
     & x & \mapsto & f(x) - \langle x_{0}^{*}, x \rangle.
\end{array}
\]
It is clear that $\phi$ is proper, and that $\text{Dom}(\phi) = \text{Dom}(f)$. Moreover $\phi$ is lower semicontinuous, since it is a composition of a lower semicontinuous function with a continuous one. We want to see now that $\phi$ is bounded below. Indeed, for all $\epsilon > 0$, since $x_{0}^{*} \in \partial_{\epsilon}f(x_{0})$ we have that
\[
\langle x_{0}^{*}, x - x_{0} \rangle \leq f(x) - f(x_{0}) + \epsilon \quad \forall \, x \in X.
\]
Therefore
\[
f(x_{0}) -  \langle x_{0}^{*}, x_{0} \rangle \leq f(x) - \langle x_{0}^{*}, x \rangle + \epsilon \quad \forall \, \epsilon > 0, \hspace{0.1cm} x \in X,
\]
which means that
\[
\phi(x_{0}) \leq \phi(x) + \epsilon \quad \forall \, \epsilon > 0, \hspace{0.1cm} x \in X.
\]
From here we deduce that $\phi$ is bounded below and in particular
\[
\phi(x_{0}) \leq \inf_{x \in X}\phi(x) + \epsilon.
\]
It is also clear that $\phi$ is $G$-invariant since it is a composition of $G$-invariant functions, and the only thing that remains to check is that $\phi$ is convex with respect to $G$. We want to see that
\[
\phi(\overline{x}) \leq \int_{G}\phi(g(x))d\mu(g),
\]
which is equivalent, by definition of $\phi$, to
\[
f\left(\int_{G}g(x)d\mu(g)\right) - \left\langle x_{0}^{*}, \int_{G}g(x)d\mu(g) \right\rangle \leq \int_{G}f(g(x)) - \langle x_{0}^{*}, g(x) \rangle d\mu(g).
\]
But this is clear, since $f$ is convex with respect to $G$ and $x_{0}^{*}$ is linear.

So, we can apply Theorem \ref{Principi variacional de Ekeland} to obtain that there exists a point $z \in \text{Dom}(\phi)$ such that
\begin{enumerate}
    \item $z$ is $G$-invariant, 
    
    \item $\phi(z) < \phi(x) + \lambda \Vert x - z \Vert$, for all $z \neq x \in X$,
    
    \item $\Vert z - x_{0} \Vert \leq \frac{\epsilon}{\lambda}$.
\end{enumerate}
By $(3)$, it is clear that $\Vert z - x_{0} \Vert \leq \frac{\epsilon}{\lambda}$, so we only need to show that $x^{*} \in \partial f(z)$, and $\Vert x^{*} - x_{0}^{*} \Vert \leq \lambda$. We define now the function
\[
\begin{array}{cccc}
    h\colon & X & \to & \mathbb{R} \\
     & x & \mapsto & \lambda\Vert x - z \Vert,
\end{array}
\]
which is proper, convex, lower semicontinuous and satisfies that $h(z) = 0$. Observe that by $(2)$ we have that
\[
(\phi + h)(x) - (\phi + h)(z) > 0 \quad \forall \, x \in X,
\]
in particular $0 \in \partial(\phi + h)(z)$. But we know that $h$ has a point of continuity in $\text{Dom}(h)\cap\text{Dom}(\phi) = \text{Dom}(\phi)$, so applying now Theorem \ref{Moreau-Rockafellar} we have that
\[
\partial(\phi + h)(z) = \partial\phi(z) + \partial h(z).
\]
Hence, $0 \in \partial\phi(z) + \partial h(z)$, therefore there exists $-z^{*} \in \partial\phi(z)$ such that $z^{*} \in \partial h(z)$. Note that the symmetrized point $\overline{z^{*}}$ belongs to $\partial h(z)$ because the subdifferential is a convex set. Moreover, we know by Corollary \ref{Corollary Moreau-Rockafellar} that $\partial\phi(z) = \partial f(z) + x_{0}^{*}$, so there exists $x^{*} \in \partial f(z)$ such that $x^{*} = -\overline{z^{*}} + x_{0}^{*}$. Notice that $x^{*}$ is $G$-invariant because so are $x_{0}^{*}$ and $\overline{z^{*}}$, and notice also that
\begin{align*}
\partial h(z) & = \left\{z^{*} \in X^{*} ~~ | ~~ \langle z^{*}, x - z \rangle \leq h(x) - h(z) \hspace{0.3cm} \forall \, x \in X\right\} \\
& = \left\{z^{*} \in X^{*} ~~ | ~~ \langle z^{*}, x - z \rangle \leq \lambda\Vert x - z \Vert \hspace{0.3cm} \forall \, x \in X\right\} \\
& = \left\{z^{*} \in X^{*} ~~ | ~~ \Vert z^{*} \Vert \leq \lambda\right\}.
\end{align*}
Therefore
\[
\Vert x^{*} - x_{0}^{*} \Vert = \Vert \overline{-z^{*}} \Vert = \Vert \overline{z^{*}} \Vert \leq \lambda.
\]
\end{proof}

\end{document}